\documentclass[reqno,a4paper]{amsart}
\usepackage{amsmath}
\usepackage{amssymb,amsfonts}
\usepackage{amsthm}
\usepackage{enumerate}
 
\usepackage[mathscr]{eucal}
\usepackage{xcolor}

\theoremstyle{plain}
\newtheorem{theorem}{Theorem}[section]

\newtheorem{conj}{Conjecture}[section]
\theoremstyle{definition}

\theoremstyle{remark}

\numberwithin{equation}{section}

\begin{document}

\title[A brief survey on the generalized Lebesgue-Ramanujan-Nagell equation]
{A brief survey on the generalized Lebesgue-Ramanujan-Nagell equation}

\author{Maohua Le and G\"{o}khan Soydan}

\address{{\bf Maohua Le}\\
Institute of Mathematics, Lingnan Normal College\\
Zhangjiang, Guangdong, 524048 China}
		
\email{lemaohua2008@163.com}

\address{{\bf G\"{o}khan Soydan} \\
	Department of Mathematics \\
	Bursa Uluda\u{g} University\\
	16059 Bursa, Turkey}
\email{gsoydan@uludag.edu.tr }

\newcommand{\acr}{\newline\indent}


\begin{abstract}
The generalized Lebesgue-Ramanujan-Nagell equation is an important type of polynomial-exponential Diophantine equation in number theory. In this survey, the recent results and some unsolved problems of this equation are given.
\end{abstract}

\maketitle

\textbf{Contents}

\bigskip

$1$ Introduction

$2$ The generalized Ramanujan-Nagell equation

$2.1$ The equation $x^2+D=k^n$

$2.2$ The equation $x^2-D=k^n$

$2.3$ The equation $x^2\pm D^m=k^n$

$3$ The generalized Lebesgue-Nagell equation

$3.1$ The equation $x^2\pm D=y^n$

$3.2$ The equation $x^2\pm D^m=y^n$

$3.3$ The divisibilty of class number of $\mathbb{Q}(\sqrt{\pm (a^2-\delta c^n)})$

References

\section{Introduction}  
A great deal of number theory arises from the discussion of the integer or rational solutions of a polynomial equation with integer coefficients. Such equations are called Diophantine equations, those Diophantine equations with variable exponents are called polynomial-exponential equations.

Let $\mathbb{Z}$, $\mathbb{N}$, $\mathbb{Q}$ denote the sets of all integers, positive integers and rational numbers, respectively. Let $D$ be a fixed positive integer, and let $p$ be a fixed prime with $\gcd(D,p)=1$. We first introduce two classical results about polynomial-exponential equations.

In 1844, E. Catalan \cite{Cat} conjectured that the equation
\begin{equation}\label{eq.1.1}
x^m+1=y^n,\, \, x,y,m,n \in \mathbb{N}, \, \, \min\{x,y,m,n\}>1 
\end{equation} 
has only the solution $(x,y,m,n)=(2,3,3,2).$\footnote{This conjecture has been proved by P. Mih\v{a}ilescu \cite{Mih}} Six years later, V.A. Lebesgue \cite{Leb} solved \eqref{eq.1.1} for $m \equiv0\pmod{2}$, namely, he proved that the equation
\begin{equation}\label{eq.1.2}
x^2+1=y^n, \, x,y,n\in\mathbb{N}, \, \min\{ x,y,n\}>1
\end{equation}
has no solutions $(x,y,n)$. Afterwards, T. Nagell \cite{Nag1,Nag3} systematically studied the solution of equations with the form
\begin{equation}\label{eq.1.3}
x^2+D=y^n, \, x,y,n\in\mathbb{N}, \, \gcd(x,y)=1,\, n>2.
\end{equation}
Therefore, \eqref{eq.1.3} and its generalizations are called the generalized Lebesgue-Nagell equation.

In 1913, S. Ramanujan \cite{Ram} conjectured that the equation
\begin{equation}\label{eq.1.4}
x^2+7=2^{n+2}, \, x,n\in\mathbb{N}
\end{equation}
has only the solutions $(x,n)=(1,1), (3,2),(5,3), (11,5)$ and $(181,13)$. In 1945, W. Ljunggren \cite{Lj3} proposed the same problem in Norwegian and T. Nagell \cite{Nag2} first proved it three years later. Therefore, the equation
\begin{equation}\label{eq.1.5}
x^2+D=
\begin{cases}
2^{n+2}, \quad  \textrm{if $p=2$},\\
\; \quad \quad \quad \quad \quad \quad \; \; \; x,n\in\mathbb{N}\\
p^{n},   \quad  \quad   \textrm{if $p\neq2$},
\end{cases}
\end{equation}
and its generalizations are called the generalized Ramanujan-Nagell equation. Moreover, since \eqref{eq.1.3} and \eqref{eq.1.5} are similar in form, they are collectively referred to as the generalized Lebesgue-Ramanujan-Nagell equation.

The methods of proving the early results of generalized Lebesgue-Ramanujan-Nagell equation depend mainly on elementary techniques from classical number theory, algebraic number theory and from Diophantine approximation. These results can be seen in the survey papers \cite{Abu8}, \cite{BP3}, \cite{Coh}, \cite{Le27}, \cite{Le44}, \cite{LeHu}, \cite{Ramas} and \cite{SarSr2}.

In this survey, we introduce some recent results concerning the generalized Lebesgue-Ramanujan-Nagell equation. In addition to the classical methods, the proof of these results used the following four tools.
\bigskip

\textbf{(i) The Baker method}
\quad

In 1966, A. Baker \cite{Bak1} first obtained an effectively computable lower bound for linear forms in the logarithms of multiple algebraic numbers. This important result and its extensions have wide application in many fields including Diophantine equations (see \cite{Bak2}, \cite{BakMas}, \cite{Bug4}, \cite{Mass}, \cite{Le44} and \cite{ShTi}).
\bigskip

\textbf{(ii) The formula for the solutions of $D_1X^2-D_2Y^2=\lambda k^z$}

Let $D_1, D_2, k$ be fixed integers such that $D_1>0$, $k>1$, $\gcd(D_1,D_2)=\gcd(D_1D_2,k)=1$ and $D_1D_2$ is not a square, and let $\lambda\in\{1,2,4\}$ or $\lambda=4$ according to $k\equiv1\pmod{2}$ or not. Around 1986, using the properties on the repsentation of integers by binary quadratic primitive forms, M.-H. Le \cite{Le3,Le11} gave the explicit formula for all solutions $(X,Y,Z)$ of the equation
\begin{equation}\label{eq.1.6}
D_1X^2-D_2Y^2=\lambda k^z, \, X,Y,Z\in\mathbb{Z},\, \gcd(X,Y)=1,\, Z>0.
\end{equation}
Because the solutions of many types of generalized Lebesgue-Ramanujan-Nagell equation are special solutions of \eqref{eq.1.6}, so this formula is pretty useful. The English version of Le's work above can be found in \cite{HL} and \cite{Le30}. In addition, Z.-F. Cao, C.I. Chu and W.C. Chiu \cite{CCC}, R.A. Mollin \cite{Mol3,Mol4}, P.-Z. Yuan \cite{Yu4, Yu7} used different methods to prove similar results. H. Yang and R.-Q. Fu \cite{YanFu1} improved the upper bound for least solutions of \eqref{eq.1.6} given in \cite{Le30}.  
\bigskip

\textbf{(iii) The BHV theorem}

Let $\alpha, \beta$ be algebraic integers. If $\alpha+\beta$ and $\alpha\beta$ are nonzero coprime integers and $\alpha/\beta$ is not a root of unity, then $(\alpha,\beta)$ is called a Lucas pair. Further, one defines the corresponding sequence of Lucas numbers by
\begin{equation}\label{eq.1.7}
L_t(\alpha,\beta)=\dfrac{\alpha^t-\beta^t}{\alpha-\beta}, \, t=0, 1,2,\cdots.
\end{equation}
A prime $q$ is called a primitive divisor of $L_t(\alpha,\beta)$ $(t>1)$ if $L_t(\alpha,\beta)\equiv0\pmod{q}$ and $(\alpha-\beta)^2L_1(\alpha,\beta)\cdots$ $L_{t-1}(\alpha,\beta)\not\equiv 0\pmod{q}.$ Similarly, we can define the Lehmer pair, the Lehmer number and its primitive divisors. In 2001, Y.F. Bilu, G. Hanrot and P.M. Voutier \cite{BHV} proved that if $t>30$, then $L_t(\alpha,\beta)$ must have primitive divisors. Thus, combine this result together with the results for $t\le 30$ due to P.M. Voutier \cite{Vout} and M. Abouzaid \cite{Abou}, the existence of primitive divisors of Lucas and Lehmer numbers is completely solved, which is usually called the BHV theorem.
\bigskip

\textbf{(iv) Modular approach}

Around 1984, G. Frey \cite{Frey} first discovered the connection between certain Diophantine equations and modularity of elliptic curves over $\mathbb{Q}$. In 1995, along this approach, A. Wiles \cite{TWill,Wiles} succesfully solved the famous Fermat's last theorem. This kind of methods are collectively referred to as modular approach. So far this method has been widely used in solving many Diophantine equations including the generalized Lebesgue-Nagell equation (see \cite{HCoh}, \cite{Falt} and \cite{Si2}).

We now illustrate the use of the above tools with the following two examples.

First, we consider \eqref{eq.1.5} for odd primes $p$. Obviously, \eqref{eq.1.5} has a solution $(x,n)$ if and only if the equation
\begin{equation}\label{eq.1.8}
X^2+DY^2=p^z, \, X,Y,Z\in\mathbb{Z},\, \gcd(X,Y)=1,\, Z>0
\end{equation}
has the solution
\begin{equation}\label{eq.1.9}
(X,Y,Z)=(x,1,n)
\end{equation}
with $Y=1$. By \cite{Le1}, if \eqref{eq.1.8} has solutions $(X,Y,Z)$, then it has a unique solution $(X_1,Y_1,Z_1)$ such that $X_1>0$, $Y_1>0$ and $Z_1\leq Z$, where $Z$ through all solutions $(X,Y,Z)$ of \eqref{eq.1.8}. The solution $(X_1,Y_1,Z_1)$ is called the least solution of \eqref{eq.1.8}. Moreover, every solution $(X,Y,Z)$ of \eqref{eq.1.8} can be expressed as
\begin{equation}\label{eq.1.10}
\begin{aligned}
&Z=Z_1t, \, t\in\mathbb{N},\\
&X+Y\sqrt{-D}=\lambda_1(X_1+\lambda_2Y_1\sqrt{-D})^t, \, \lambda_1, \lambda_2 \in\{1,-1\}.
\end{aligned}
\end{equation} 
Therefore, applying \eqref{eq.1.10} to \eqref{eq.1.9}, we find that \eqref{eq.1.5} has solutions $(x,n)$ if and only if the least solution $(X_1,Y_1,Z_1)$ satisfies $Y=1$. Obviously, if $Y_1=1$, then \eqref{eq.1.5} has solution $(x,n)=(X_1,Z_1)$. Further, if \eqref{eq.1.5} has another solution $(x,n)$, then there exists an odd positive integer $t$ such that $t>1$ and
\begin{equation}\label{eq.1.11}
n=Z_1t, \,
X+\sqrt{-D}=\lambda_1(X_1+\lambda_2\sqrt{-D})^t, \, \lambda_1, \lambda_2 \in\{1,-1\}.
\end{equation} 

When $t=3$, we get from \eqref{eq.1.11} that
\begin{equation}\label{eq.1.12}
(D,p)=
\begin{cases}
	(2,3),  \\
	(3s^2+1,4s^2+1),
\end{cases}
(x,n)=
\begin{cases}
	(1,1) \, \textrm{and} \, (5,3),  \\
(s,1) \, \textrm{and} \, (8s^3+3s,3), \quad s\in\mathbb{N}.
\end{cases}
\end{equation}
Since, by \cite{Ape}, \eqref{eq.1.5} has at most two solutions $(x,n)$ for odd primes $p.$ Hence, if $D,p$ satisfy \eqref{eq.1.12}, then \eqref{eq.1.5} is solved.

When $t>3$, let
\begin{equation}\label{eq.1.13}
\alpha=X_1+\sqrt{-D}, \quad \beta=X_1-\sqrt{-D}. 
\end{equation}
Then, by \eqref{eq.1.11}, we have
\begin{equation}\label{eq.1.14}
  \left\lvert \dfrac{\alpha^t-\beta^t}{\alpha-\beta}\right\rvert=1
\end{equation}
Let $C_j$ $(j=1,2,\cdots)$ denote effectively computable absolute constants. Using the Baker method, by \eqref{eq.1.14}, we have $t<C_1$ and $\max\{D,p\}<C_2$. This implies that if $\max\{D,p\}>C_2$ and $(D,p)\neq(3s^2+1,4s^2+1)$, then \eqref{eq.1.5} has at most one solution $(x,n)$ (see \cite{Le2}, \cite{Le5} and \cite{Le21}).

In addition, \eqref{eq.1.14} gives the connection between the number of solutions of generalized Ramanujan-Nagell equation and the multiplicity of second order recurrences. The work on the latter can be seen in \cite{Beu1}, \cite{Le16}, \cite{Le50} and \cite{Yu5}.

On the other hand, we see from \eqref{eq.1.13} that $(\alpha,\beta)$ is a Lucas pair. Further, by \eqref{eq.1.7} and \eqref{eq.1.14}, the Lucas number $L_t(\alpha,\beta)$ has no primitive divisors. Therefore, using the BHV theorem, \eqref{eq.1.14} is false if $t>3$. Thus, for odd primes $p$, \eqref{eq.1.5} has at most one solution $(x,n)$, except for $D,p$ satisfy \eqref{eq.1.12}. We can get similar results for a lot of other types of generalized Ramanujan-Nagell equation in this approach (see \cite{BugSh}, \cite{Le17}, \cite{Le20}, \cite{Le29} and \cite{Le39}).

Finally, we consider \eqref{eq.1.3} for $D\equiv 7\pmod{8}$. In this case, \eqref{eq.1.3} might have solutions $(x,y,n)$ with $y\equiv 0\pmod{2}$ and $y$ is not a power of 2. This implies that the approach outlined above will not work, and the solution of \eqref{eq.1.3} become quite difficult. For example, if $D=7$, then \eqref{eq.1.3} can be written as
\begin{equation}\label{eq.1.15}
x^2+7=y^n, \, x,y,n\in\mathbb{N},\, \gcd(x,y)=1,\, n>2.
\end{equation}
Obviously, \eqref{eq.1.15} has no solutions $(x,y,n)$ with $y\equiv 1\pmod{2}$, and it has only solutions $(x,y,n)=(1,2,3), (3,2,4),(5,2,5), (11,2,7),(181,2,15),(181,8,5)$ and $(181,32,3)$ with $y$ is a power of 2 (see \cite{Le38}). In 1993, J.H.E. Cohn \cite{Cohn5} conjectured that \eqref{eq.1.15} has no other solutions $(x,y,n)$. This problem has not been solved for years. In 1998, using the Baker method, J.-L. Lesage \cite{Les} proved that all solutions $(x,y,n)$ of \eqref{eq.1.15} satisfy $n<6.6\times10^{15}$. Five years later, S. Siksek and J.E. Cremona \cite{SiCr} proved that \eqref{eq.1.15} has no solutions $(x,y,n)$ with $11\le n\le 10^8.$ It was not until 2006 that Y. Bugeaud, M. Mignotte and S. Siksek \cite{BMS2} used the modular approach to confirmed Cohn's conjecture. 
 
\section{The generalized Ramanujan-Nagell equation}

\subsection{The equation $x^2+D=k^n$}

\quad

\quad

Let $D_1$, $D_2$ be fixed positive integers with $\gcd(D_1,D_2)=1$, and let $p$ be a fixed prime with $\gcd(D_1D_2,p)=1$. In this subsection we are going to ignore the classical results of \eqref{eq.1.5} and go consider to its generalization with the following form
\begin{equation}\label{eq.2.1}
D_1x^2+D_2=
\begin{cases}
2^{n+2}, \quad  \textrm{if $p=2$},\\
\; \quad \quad \quad \quad \quad \quad \; \; \; x,n\in\mathbb{N}.\\
p^{n},   \quad  \quad   \textrm{if $p\neq2$},
\end{cases}
\end{equation}
In 1918, using Thue's method on Diophantine approximation, G. P\'{o}lya \cite{Pol} proved that $P[D_1x^2+D_2]\rightarrow \infty$ as $|x|\rightarrow \infty$, where $P[D_1x^2+D_2]$ is the great prime divisor of $D_1x^2+D_2$. This implies that \eqref{eq.2.1} has only finitely many solutions $(x,n)$. Let $N(D_1,D_2,p)$ denote the number of solutions $(x,n)$ of \eqref{eq.2.1}. The study of \eqref{eq.2.1} mainly discusses the upper bound of $N(D_1,D_2,p)$.

In 2001, using the BHV theorem, Y. Bugeaud and T.N. Shorey \cite{BugSh} completely determined the upper bound of $N(D_1,D_2,p)$ as follows. 
\begin{theorem}
$N(D_1,D_2,2)\le 1$, except for the following cases:

$(i)$ $N(1,7,2)=5$, $(x,n)=(1,1),(3,2),(5,3),(11,5)$ and $(181,13).$

$(ii)$ $N(3,5,2)=3$, $(x,n)=(1,1),(3,3)$ and $(13,7).$

$(iii)$ $N(1,23,2)=2$, $(x,n)=(3,3)$ and $(45,9).$

$(iv)$ $N(1,2^{r+2}-1,2)=2$, $(x,n)=(1,r)$ and $(2^{r+1}-1,2r),$
where $r$ is a positive integer with $r\ge 2$.

$(v)$ $N(3,29,2)=2$, $(x,n)=(1,3)$ and $(209,15)$.

$(vi)$ $N(13,3,2)=2$, $(x,n)=(1,2)$ and $(71,14)$.

$(vii)$ $N(21,11,2)=2$, $(x,n)=(1,3)$ and $(79,15)$.

$(viii)$ If $D_1s^2=2^r-\delta$, $D_2=3\cdot2^r+\delta$, where $r,s$ are positive integers with $r\ge 2$, $\delta\in\{1,-1\}$, then $N(D_1,D_2,2)=2$, $(x,n)=(s,r)$ and $((2^{r+1}+\delta)s,3r)$.
\end{theorem}

\begin{theorem}
For odd primes, if $4D_1s^2=p^r-\delta$, $4D_2=3p^r+\delta$, where $r,s$ are positive integers, $\delta\in\{1,-1\}$, then $N(D_1,D_2,p)=2$, $(x,n)=(s,r)$ and $((2p^{r}+\delta)s,3r)$. Otherwise, $N(D_1,D_2,p)\le 1$.
\end{theorem}

For any nonnegative integer $i$, let $F_i$ and $L_i$ denote the $i-th$ original Fibonacci and Lucas numbers, rspectively. We know introduce their two properties as follows:

$(i)$ (J.H.E. Cohn \cite{Cohn1}) The equation
\begin{equation*}
L_i=2z^2,\, z\in\mathbb{N}
\end{equation*}
has only the solution $(i,z)=(0,1)$. 

$(ii)$ (Y. Bugeaud, M. Mignotte and S. Siksek \cite{BMS1}) The equation
\begin{equation*}
F_i=z^r,\, z,r\in\mathbb{N}, \, z>1, r>1
\end{equation*}
has only the solutions $(i,z,r)=(6,2,3)$ and $(12,12,2)$.

For odd positive integers $D_1$, $D_2$ and odd primes $p$, let $N(D_1,D_2,2,p)$ and $N(D_1,D_2,4,p)$ denote the number of solutions $(x,n)$ of the equations
\begin{equation}\label{eq.2.2}
D_1x^2+D_2=2p^n, \, x,n\in\mathbb{N}
\end{equation}
and 
\begin{equation}\label{eq.2.3}
D_1x^2+D_2=4p^n, \, x,n\in\mathbb{N},
\end{equation}
respectively. By the results of \cite{BugSh} with the above mentioned properties of Fibonacci and Lucas numbers, we have

\begin{theorem}
$N(D_1,D_2,2,p)\le 1$, except for the following cases:

$(i)$ $N(1,1,2,5)=2$, $(x,n)=(3,1)$ and $(7,2).$

$(ii)$ $N(1,1,2,13)=2$, $(x,n)=(5,1)$ and $(239,4).$

$(iii)$ $N(7,11,2,3)=2$, $(x,n)=(1,2)$ and $(1169,14).$

$(iv)$ $N((p^r-\delta)/2s^2,(3p^r+\delta)/2,2,p)=2$, $(x,n)=(s,r)$ and $((2p^r+\delta)s,3r),$ where $r,s$ are positive integers, $\delta\in\{1,-1\}$.

\end{theorem}

\begin{theorem}
$N(D_1,D_2,4,p)\le 1$, except for the following cases:
	
$(i)$ $N(1,11,4,3)=3$, $(x,n)=(1,1),(5,2)$ and $(31,5)$.
	
$(ii)$ $N(1,19,4,5)=3$, $(x,n)=(1,1),(9,2)$ and $(559,5)$
	
$(iii)$ $N(1,3,4,7)=2$, $(x,n)=(5,1)$ and $(37,3).$
	
$(iv)$ $N(1,4p^r-1,4,p)=2$, $(x,n)=(1,r)$ and $(2p^r-1,2r)$, where $r$ is a positive integer.

$(v)$ $N(F_{6t-3}/s^2,L_{6t},4,F_{6t-1})=N(F_{6t-1}/s^2,L_{6t+2},4,F_{6t+1})=2$, $(x,n)=(s,1)$ and $(s(D_1^2s^4-5D_1s^2p+5p^2),5)$, where $t$ is a positive integer.
\end{theorem}

Let $k$ be a fixed positive integer such that $k>1$ and $\gcd(D_1D_2,k)=1$. For any positive integer $a$ with $a>1$, let $\omega(a)$ denote the number of distinct prime divisors of $a$. Then the equation
\begin{equation}\label{eq.2.4}
D_1x^2+D_2=
\begin{cases}
\lambda k^{n}, \quad  \textrm{if $k\equiv 1\pmod{2}$},\quad  \lambda\in\{1,2,4\},\\
\; \quad \quad \quad \quad \quad \quad \quad \quad \quad \quad \quad \quad \quad \quad \quad \quad \quad \quad \quad x,n\in\mathbb{N}\\
4k^{n},   \quad    \textrm{if $k\equiv 0\pmod{2}$},
\end{cases}
\end{equation}
is a naturally generalization of \eqref{eq.2.1}, \eqref{eq.2.2} and \eqref{eq.2.3}. Obviously, if \eqref{eq.2.4} has a solution $(x,n)$, then the equation 
\begin{equation}\label{eq.2.5}
D_1x^2+D_2Y^2=
\begin{cases}
\lambda k^{z},\\ 
\quad \quad  X,Y,Z\in\mathbb{Z}, \, \gcd(X,Y)=1, Z>0\\
4k^{z},
\end{cases}
\end{equation}
has a corresponding solution $(X,Y,Z)=(x,1,n)$ with $Y=1$. By \cite{Le30}, all solutions $(X,Y,Z)$ of \eqref{eq.2.5} can be put into at most $2^{\omega(k)-1}$ classes, and every class, say $S$, has a unique solution $(X_1,Y_1,Z_1)$ such that $X_1>0$, $Y_1>0$ and $Z_1\le Z$, where $Z$ through all solutions $(X,Y,Z)$ in $S$. The solution $(X_1,Y_1,Z_1)$ is called the least solution of $S$. It can be seen from the analysis in Introduction that if \eqref{eq.2.4} has a solution $(x,n)$ which make the corresponding solution $(X,Y,Z)=(x,1,n)$ of \eqref{eq.2.5} belongs to $S$, then the least solution $(X_1,Y_1,Z_1)$ of $S$ satisfies $Y_1=1$. Let $N(D_1,D_2,k)$ denote the number of solutions $(x,n)$ of \eqref{eq.2.4}. Let $C_j(a,b,\cdots)$ $(j=1,2,\cdots)$ denote effectively computable constants depending only on the parameters $a,b,\cdots$. Using the Baker method, E.A. Bender and N.P. Herzberg \cite{BH} proved that if $k\equiv 1 \pmod{2}$ and $k>C_1(D_1D_2)$, then $N(D_1,D_2,k)\le 2^{\omega(k)-1}$. T.-J. Xu and M.-H. Le \cite{XuLe} showed that, except for some explicit exceptional cases, the above condition $k>C_1(D_1D_2)$ can be improved to $\max\{D_1D_2,k\}>C_3.$ Obviously, using the BHV theorem, these conditions can be removed (see \cite{BugSh}). In addition, X.-G. Chen, Y.-D. Guo and M.-H. Le \cite{CGL} proved that if $D_1=1$, $D_2$ is an odd prime, $k\equiv 1 \pmod{2}$ and $\lambda=1$, then $N(D_1,D_2,k)\le 8.$

As a special case of \eqref{eq.2.4}, Z.-G. Li and P.-Z. Yuan \cite{LiYu} claimed that if $\lambda=4$ and $(D_1,D_2)=(k-1,91k+9)$, then \eqref{eq.2.4} has only the solution $(k,x,n)=(5,3,3)$. However, M.-H. Le \cite{Le58} pointed that the solutions $(k,x,n)=(7,11,3)$ and $(25,3,4)$ were left out in \cite{LiYu}, and proved some further results.

The equations \eqref{eq.2.1}-\eqref{eq.2.4} relate to a classical problem concerning the equation
\begin{equation}\label{eq.2.6}
\dfrac{X^m-1}{X-1}=\dfrac{Y^n-1}{Y-1},\, X,Y,m,n\in\mathbb{N},\, X>Y>1,\, n>m>2.
\end{equation}
Around 1916, R. Ratat \cite{Rat} and Goormaghtigh \cite{Goor} conjectured that \eqref{eq.2.6} has only the solutions $(X,Y,m,n)=(5,2,3,5)$ and $(90,2,3,13)$. In addition, P.T. Bateman asked that if $(X,Y,m,n)=(5,2,3,5)$ the only solution of \eqref{eq.2.6} for which $X,Y$ and $(X^m-1)/(X-1)$ are all primes? (see Problem B 25 of \cite{Guy}). These are two problems have not solved as yet. Since the known solutions of \eqref{eq.2.6} satisfy $m=3$, it is a most interesting case of \eqref{eq.2.6}. Then, \eqref{eq.2.6} can be written as
\begin{equation}\label{eq.2.7}
(Y-1)(2X+1)^2+(3Y+1)=4Y^n,\, X,Y,n\in\mathbb{N}, \, X>Y>1,\, n>3.
\end{equation}  
In this respect, there have the following results:

$(i)$ (M.-H. Le \cite{Le45}) \eqref{eq.2.7} has only the solutions $(X,Y,n)=(5,2,5)$ and $(90,2,13)$ with $Y$ is a prime power.

$(ii)$ (P.-Z. Yuan \cite{Yu8}; B. He \cite{Bhe}) \eqref{eq.2.7} has only the solutions $(X,Y,n)=(5,2,5)$ and $(90,2,13)$ with $n\equiv 1 \pmod{2}$.

$(iii)$ (M.-H. Le \cite{Le49}) If $(X,Y,n)$ is a solution of \eqref{eq.2.7} with $(X,Y,n)\neq (5,2,5)$ or $(90,2,13)$, then $\gcd(X,Y)>1$ and $X\not\equiv 0 \pmod{Y}$.

Let $D,k$ be fixed odd positive integers with $\gcd(D,k)=1$. The equation
\begin{equation}\label{eq.2.8}
x^2+D=k\cdot 2^n,\, x,n\in\mathbb{N}
\end{equation}
is another generalization of \eqref{eq.1.5}. In this respect, M.A. Bennett, M. Filaseta and O. Trifonov \cite{BFT} proved that if $D=7$, then the solutions $(x,n)$ of \eqref{eq.2.8} satisfy $\sqrt{x}<k$. J. Stiller \cite{Stil} proved that if $(D,k)=(119,15)$, then \eqref{eq.2.8} has six solutions $(x,n)=(1,3).(11,4),(19,5),(29,6),(61,8)$ and $(701,15)$. Afterwards, M. Ulas \cite{Ul} interested in finding equations of the form \eqref{eq.2.8} with many solutions. For example, he showed that if $(D,k)=(2^{3r+3}-1,2^{3r}+1)$, where $r$ is a positive integer, then \eqref{eq.2.8} has at least five solutions $(x,n)=(3,3),(2^{2r+1}-2^{r+1}-1,r+2)$, $(2^{3r+1}-1,3r+2)$, $(2^{3r+2}+1,3r+4)$ and $(2^{6r+3}+2^{3r+2}-1,9r+6)$. In the same paper, he put forward a lot of conjectures in this respect. Some of them have been confirmed by Z.-F. Zhang and A. Togb\'{e} \cite{ZT}.

Let $p_1,p_2,\cdots,p_r$ $(r>1)$ be distinct primes with $p_i\nmid D$ $(i=1,2,\cdots,r)$. By the result of G. P\'{o}lya \cite{Pol}, we can find that the equation
\begin{equation*}
x^2+D=p_1^{n_1}p_2^{n_2}\cdots p_r^{n_r}, \, x\in \mathbb{N},\, n_1,n_2,\cdots n_r\in \mathbb{Z},\, n_i\ge 0,\, i=1,2,\cdots,r 
\end{equation*}
has only finitely many solutions $(x,n_1,n_2,\cdots,n_r)$. For example T. Yamada \cite{Yamada} proved that if $p_1=2$ and $r=3$, then the equation has at most 63 solutions $(x,n_1,n_2,n_3)$ with $0 \le n_1 \le 2$.
\subsection{The equation $x^2-D=k^n$}
\quad

\quad

Let $D$ be fixed positive integer integer which is not a square, and let $p$ be a fixed prime with $\gcd(D,p)=1$. In this subsection we consider the equation
\begin{equation}\label{eq.2.9}
x^2-D=
\begin{cases}
2^{n+2}, \quad  \textrm{if $p=2$},\\
\; \quad \quad \quad \quad \quad \quad \; \; \; x,n\in\mathbb{N}\\
p^{n},   \quad  \quad   \textrm{if $p\neq2$},
\end{cases}
\end{equation} 
and its generalizations. Obviously, according to the language of algebraic number theory, the left side of \eqref{eq.1.5} can be splitten in imaginary quadratic fields, and the left side of \eqref{eq.2.9} must be splitten in real quadratic fields. This difference makes the solution of \eqref{eq.2.9} much more difficult than \eqref{eq.1.5}.

Let $N(-D,p)$ denote the number of solutions $(x,n)$ of \eqref{eq.2.9}. In 1980, using the Diophantine approximation method based on hypergeometric series, F. Beukers \cite{Beu2} proved that $N(-D,2)\le 4$. In the same paper, he showed that if $N(-D,2)>3$, the $D$ must be one of the following three cases:

$(i)$ If $D=2^{2r}-3\cdot2^{r+1}+1$, where $r$ is a positive integer with $r\ge 3$, then $N(-D,2)=4$, $(x,n)=(2^r-3,1)$, $(2^r-1,r)$, $(2^r+1,r+1)$ and $(3\cdot2^r-1,2r+1)$

$(ii)$ If $D=(\frac{1}{3}(2^{2r+1}-17))^2-32$, where $r$ is a positive integer with $r\ge 3$, then \eqref{eq.2.9} has at least three solutions $(x,n)=(\dfrac{1}{3}(2^{r+1}-17),3)$, $(\dfrac{1}{3}(2^{r+1}+1),2r+1)$ and $(\dfrac{1}{3}(17\cdot2^{r+1}-1),4r+5).$

$(iii)$ If $D=2^{2r_1}+2^{2r_2}-2^{r_1+r_2+1}-2^{r_1+1}-2^{r_2+1}+1$, where $r_1$, $r_2$ are positive integers with $r_2>r_1+1>2$, then \eqref{eq.2.9} has at least three solutions $(x,n)=(2^{r_2}-2^{r_1}-1,r_1)$, $(2^{r_2}-2^{r_1}+1,r_2)$ and $(2^{r_2}+2^{r_1}-1,r_1+r_2)$.
			
In 1992, M.-H. Le \cite{Le19} further proved the following result:

\begin{theorem}
$N(-D,2)\le 3$, except for $D$ satisfies the above condition $(i)$.
\end{theorem}

In addition, F. Beukers put forward the following conjecture in his letter to M.-H. Le.
\begin{conj}\label{conj.2.1}
$N(-D,2)\le 2$, except for $D$ satisfies the above conditions $(i)$, $(ii)$ and $(iii).$
\end{conj}

In this respect, M.-H. Le \cite{Le14} proved that if the equation
\begin{equation}\label{eq.2.10}
U^2-DV^2=-1,\, U,V\in\mathbb{N}
\end{equation} 
has solutions $(U,V)$, then Conjecture \ref{conj.2.1} is true. S.-C. Yang \cite{SYang2}, H.-Y. Chen \cite{HCh} discussed the upper bound of $N(-D,2)$ for the case that $D=2^rs+1$ where $r$ and $s$ are positive integers with $s\equiv 1 \pmod{2}$. But, in general, Conjecture \ref{conj.2.1} is still an unsolved problem.

For odd primes $p$, using the same method as in \cite{Beu2}, F. Beukers \cite{Beu3} proved that $N(-D,p) \le 4$. Simultaneously, he proposed the following conjecture:
\begin{conj}\label{conj.2.2}
For odd primes $p$, $N(-D,p)\le 3$.
\end{conj}

In 1991, using the Baker method, M.-H. Le \cite{Le13} basically confirmed Conjecture \ref{conj.2.2} that if $max\{D,p\}>10^{240}$, then $N(-D,p)\le 3$. Afterwards, the condition $\max\{D,p\}>10^{240}$ has been continously improved by M.-H. Le \cite{Le15,Le26} and P.-Z. Yuan \cite{Yu3}. Until 2002, using the Diophantine approxmation method, M. Bauer and M.A. Bennett \cite{BB1} completely confirmed Conjecture \ref{conj.2.2}, namely, they proved the following result:
\begin{theorem}
For odd primes $p$, $N(-D,p)\le 3$.
\end{theorem}

Because  we know that $N(-D,p)\le 3$, by \cite{Beu3}, we have $N(-D,p)=3$ if $D,p$ satisfy one of the following two cases:

$(i)$ $(D,p)=(\dfrac{1}{4}(3^{2r+1}+1))^2-3^{2r+1},3)$, $(x,n)=(\dfrac{1}{4}(3^{2r+1}-7),1),$ $(\dfrac{1}{4}(3^{2r+1}+1),2r+1)$ and $(2\cdot 3^{2r+1}-\dfrac{1}{4}(3^{2r+1}+1),4r+3)$, where $r$ is a positive integer.

$(ii)$ $(D,p)=((\dfrac{1}{4s}(p^{r}-1))^2-p^r,4s^2+1)$, $(x,n)=(\dfrac{1}{4s}(p^{r}-1)-2s,1)$, $(\dfrac{1}{4s}(p^{r}-1),r)$ and $(2sp^r+\dfrac{1}{4s}(p^{r}-1),2r+1)$ where $r$, $s$ are positive integers.

Since no other pairs $(D,p)$ have been found that would make $N(-D,p)=3$, M.-H. Le \cite{Le44} proposed a further conjecture as follows.

\begin{conj}\label{conj.2.3}
For odd primes $p$, $N(-D,p)\le 2$, except for $D$, $p$ satisfy the above cases $(i)$ and $(ii)$.
\end{conj}

In this respect, T.-T. Wang and Y.-Z. Jiang \cite{WangJia} proved that $N(-2^r,p)\le 1$, where $r$ is a positive integer. F. Luca \cite{Luc3} proved that $N(-(p^r+1),p)=0$, where $r$ is a positive integer. J.-M. Yang \cite{JYang} proved that if $D>10^{12}$ and \eqref{eq.2.10} has solutions $(U,V)$, then $N(-D,3)\le 2$. Y.-Y. Qu, H. Gao and Q.-W. Mu \cite{QCM} improved the condition $D>10^{12}$ to $D>10^6$. Y.-E. Zhao and T.-T. Wang \cite{ZhWang} further proved that if $D>C_2(p)$, \eqref{eq.2.10} has solutions $(U,V)$ and the least solution $(u_1,v_1)$ of Pell's equation
\begin{equation}\label{eq.2.11}
u^2-Dv^2=1, \, u,v\in\mathbb{N}
\end{equation}
stisfies $v_1\not\equiv 0 \pmod{p}$, then $N(-D,p)\le 2$. But, in general, this problem has not yet been solved.

For odd positive integers $D$ and odd primes $p$, we now consider the equation
\begin{equation}\label{eq.2.12}
x^2-D=4p^n,\, x,n\in\mathbb{N} 
\end{equation}
Let $N(-D,4,p)$ denote the number of solutions $(x,n)$ of \eqref{eq.2.12}. We first deal with the case that $D=4p^r+1$, where $r$ is a positive integer. Obviously, this case is related to the equation
\begin{equation}\label{eq.2.13}
x^2=4q^m+4q^l+1, \, x,m,l\in\mathbb{N},\, m>l,\, \gcd(m,l)=1,
\end{equation} 
where $q=p^d$, $d$ is a positive integer. In this respect, A. Bremner, R. Calderbank, P. Hanlon, P. Morton and J. Wolfskill \cite{BCHMW} proved that if $q=3$ then \eqref{eq.2.13} has only the solution $(x,m,l)=(5,1,1),(7,2,1)$ and $(11,3,1)$ with $l=1$. R. Calderbank \cite{Cal} conjectured that if $q\neq 3$, then \eqref{eq.2.13} has only the solution $(x,m,l)=(2q+1,2,1)$ with $l=1$. N. Tzanakis and J. Wolfskill \cite{TzWol1,TzWol2} confirmed this conjecture, and proved that \eqref{eq.2.13} has no solution $(x,m,l)$ with $l=2$. Afterwards, M.-H. Le \cite{Le9} further proved that \eqref{eq.2.13} has no solutions $(x,m,l)$ with $l>2$. Thus, the above results show that $N(-(4p^r+1),4,p)=1$ and $(x,n)=(2p^r+1,2r)$, except for $N(-13,4,3)=3$, $(x,n)=(5,1),(7,2)$ and $(11,3)$.

Beside the case $D=4p^r+1$, M.-H. Le \cite{Le18} also noticed the following two cases:

$(i)$ If $D=(\dfrac{1}{2}(3^{2r}+1))^2-4\cdot3^{2r}$ and $p=3$, where $r$ is a positive integer with $r\ge 2$, then \eqref{eq.2.12} has at least three solutions $(x,n)=(\dfrac{1}{2}(3^{2r}-7),1)$, $(\dfrac{1}{2}(3^{2r}+1),2r)$ and $(\dfrac{1}{2}(7\cdot 3^{2r}-1),4r+1)$.

$(ii)$ If $D=p^{2r_2}+p^{2r_1}-2p^{r_1+r_2}-2p^{r_1}-2p^{r_2}+1$, where $r_1$, $r_2$ are positive integers with $r_2>r_1$, then \eqref{eq.2.12} has at least three solutions $(x,n)=(p^{r_2}-p^{r_1}-1,r_1)$, $(p^{r_2}-p^{r_1}+1,r_2)$ and  $(p^{r_2}+p^{r_1}-1,r_1+r_2).$

In \cite{Le18}, using the Baker method, he proved that
\begin{equation*}
N(-D,4,p)\le
\begin{cases}
5, \quad  \textrm{if $(D,p)$ is of the case $(i)$}, \\
4,   \quad \quad   \textrm{otherwise}.
\end{cases}
\end{equation*}
Moreover, if $\max\{D,p\}>10^{60}$, then
\begin{equation*}
N(-D,4,p)\le
\begin{cases}
4, \quad  \textrm{if $(D,p)$ is of the case $(i)$ or $(ii)$}, \\
3,   \quad \quad   \textrm{otherwise},
\end{cases}
\end{equation*}
Simultaneously, he proposed the following conjeture:
\begin{conj}\label{conj.2.4}
For $D\neq 4p^r+1$, if $(D,p)$ is of the case $(i)$ or $(ii)$, then $N(-D,4,p)=3$. Otherwise, $N(-D,4,p)\le 2$.
\end{conj}
Let $D,k$ be fixed positive integers such that $D$ is not a square, $k>1$ and $\gcd(2D,k)=1$. Then the equation
\begin{equation}\label{eq.2.14}
x^2-D=k^n, \, x,n\in\mathbb{N}
\end{equation}
is a naturally generalization of \eqref{eq.2.9}. Let $N(-D,k)$ denote the number of solutions $(x,n)$ of \eqref{eq.2.14}. In this respect, we have the following results:

$(i)$ (X.-G. Chen and M.-H. Le \cite{CL1}) If $\max\{D,k\}>C_4$ then $N(-D,k)\le 3\cdot 2^{\omega(k)-1}+1.$

$(ii)$ (M.-H. Le \cite{Le33}) $N(-D,k)\le 2^{\omega(D)+1}+1$.

$(iii)$ (B. He and A. Togb\'{e} \cite{HT}) $N(-D,k)\le 6(\log[3.2D])/\log k+8$, where $[3.2D]$ is the integer part of $3.2D.$ 

\subsection{The equation $x^2\pm D^m=k^n$}
\quad

\quad

In this subsection we first consider an exponential generalization of \eqref{eq.1.5} with the form
\begin{equation}\label{eq.2.15}
x^2+D^m=
\begin{cases}
2^{n+2}, \quad  \textrm{if $p=2$},\\
\; \quad \quad \quad \quad \quad \quad \; \; \; x,m,n\in\mathbb{N}.\\
p^{n},   \quad  \quad   \textrm{if $p\neq 2$},
\end{cases}
\end{equation}
Let $N^*(D,p)$ denote the number of solutions $(x,m,n)$ of \eqref{eq.2.15}. In 2001, Y. Bugeaud \cite{Bug3} showed that if $D\neq 7$ or 15, then $N^*(D,2)\le 1$. Obviously, this result missed many things. Using the BHV theorem, by combining the results of \cite{Bug3}, \cite{BugSh} and \cite{Le25}, we get a complete result as follows.
\begin{theorem}
$N^*(D,2)\le 1$, except for the following cases:

$(i)$ $N^*(7,2)=6$, $(x,m,n)=(1,1,1), (3,1,2), (5,1,3), (11,1,5), (181,1,13)$ and $(13,3,7).$

$(ii)$ $N^*(23,2)=2$, $(x,m,n)=(3,1,3)$ and $(45,1,9).$

$(iii)$ $N^*(2^{r+2}-1,2)=2$, $(x,m,n)=(1,1,r)$ and $(2^{r+1}-1,1,2r)$, where $r$ is a positive integer with $r\ge 2$.
\end{theorem}

For fixed odd positive integers $D_1$, $D_2$ such that $\min\{D_1,D_2\}>1$ and $\gcd(D_1,D_2)=1$, let $N^*(D_1,D_2,2)$ denote the number of solutions $(x,m,n)$ of the equation
\begin{equation}\label{eq.2.16}
D_1x^2+D_2^m=2^{n+2}, \, x,m,n\in\mathbb{N}.
\end{equation}
Very recently, J.-H. Li \cite{Li} proved a general result as follows.
\begin{theorem}
$N^*(D_1,D_2,2)\le 2$, except for the following cases:

$(i)$ $N^*(3,5,2)=4$, $(x,m,n)=(1,1,1)$, $(3,1,3), (13,1,7)$ and $(1,3,5)$.

$(ii)$ $N^*(5,3,2)=4$, $(x,m,n)=(1,1,1)$, $(5,1,5), (1,3,3)$ and $(19,5,9)$.

$(iii)$ $N^*(13,3,2)=3$, $(x,m,n)=(1,1,2)$, $(71,1,14)$ and $(1,5,6)$.

$(iv)$ $N^*(31,97,2)=3$, $(x,m,n)=(1,1,5)$, $(65,1,15)$ and $(15,2,12)$.
\end{theorem}

For odd primes $p$, since the eighties, M. Kutsuna \cite{Kut}, Q. Sun and Z.-F. Cao \cite{SuCa}, M. Toyoizumi \cite{Toyo}, M. Yamabe \cite{Yama1, Yama2, Yama3} solved \eqref{eq.2.15} for some special cases. In general, by \cite{Bug3} and \cite{Le7}, we have $N^*(D,p)\le 2$, except for $N^*(2,3)=4$ and $D,p$ satisfy
\begin{equation}\label{eq.2.17}
D=3s^2+1, \, p=4s^2+1, \, s\in\mathbb{N}.
\end{equation}
In 2005, P.-Z. Yuan and Y.-Z. Hu \cite{YuHu1} used the BHV theorem to solve the case \eqref{eq.2.17}. They proved that $N^*(3s^2+1,4s^2+1)=2$ except for $N^*(4,5)=3$. Afterwards, using some different methods, K.-Y. Chen \cite{KCh1}, M.-H. Le \cite{Le52, Le55, Le56}, Z.-W. Liu \cite{ZLiu1}, S.-C. Yang \cite{SYang3} independently proved the same result. Thus, combine the above results, we have
\begin{theorem}
For odd primes, $N^*(D,p)\le 2$, except for the following cases:

$(i)$ $N^*(2,3)=4$, $(x,m,n)=(1,1,1)$, $(5,1,3), (1,3,2)$ and $(1,11,5)$.

$(ii)$ $N^*(2,5)=3$, $(x,m,n)=(1,2,1)$, $(11,2,3)$ and $(3,4,2)$.

$(iii)$ $N^*(4,5)=3$, $(x,m,n)=(1,1,1)$, $(11,1,3)$ and $(3,2,2)$.
\end{theorem}

In addition, Y.-Z. Hu \cite{Hu2}, T. Miyazaki \cite{Miy} discussed some properties of \eqref{eq.2.15}.

For odd primes $p$, let $N^*(D_1,D_2,p)$ denote the number of solutions $(x,m,n)$ of the equation
\begin{equation}\label{eq.2.18}
D_1x^2+D_2^m=p^n, \, x,m,n\in\mathbb{N},
\end{equation}  
where $D_1$, $D_2$ are fixed positive integers such that $\min\{D_1,D_2\}>1$ and
$\gcd(D_1,D_2)=\gcd(D_1D_2,p)=1$. In 2010, P.-Z. Yuan and Y.-Z. Hu \cite{YuHu2} proved that if $D_1=\dfrac{1}{4}(p^r-1)$ and $D_2=\dfrac{1}{4}(3p^r+1)$, where $r$ is a positive integer, then $N^*(D_1,D_2,p)\le 2$, except for $N^*(2,7,3)=3$. Two years later, Y.-Z. Hu and M.-H. Le \cite{HuLe2} proved a general result as follows.
\begin{theorem}
$N^*(D_1,D_2,p)\le 2$, except for the following cases:

$(i)$ $N^*(2,7,3)=3$, $(x,m,n)=(1,1,2)$, $(19,1,6)$ and $(4,2,6)$.

$(ii)$ $N^*(10,3,13)=3$, $(x,m,n)=(1,1,1)$, $(4,2,2)$ and $(1,7,3)$.

$(iii)$ $N^*(10,3,37)=3$, $(x,m,n)=(1,3,1)$, $(71,5,3)$ and $(8,6,2)$.

$(iv)$ $N^*((3^{2r}-1)/s^2,3,4\cdot 3^{2r-1}-1)=3$, $(x,m,n)=(2s,1,1)$, $(s,2r,1)$ and $((8\cdot 3^{2r-1}+1)s,2r+2,3)$, where $r$, $s$ are positive integers.
\end{theorem}

In 1962, J.-R. Chen \cite{JCh} discussed the solutions $(x,m,n)$ of the equation
\begin{equation}\label{eq.2.19}
x^2=D^{2m}-D^mp^n+p^{2n}, \, x,m,n\in\mathbb{N}
\end{equation}
in the process of studying J\'{e}smanowicz' conjecture concerning Pythagorean triples (see \cite{Jes}). Obviously, \eqref{eq.2.19} looks like \eqref{eq.2.15}. In this respect, R.-Z. Tong \cite{Tong1} found all solutions $(x,m,n)$ of \eqref{eq.2.19} with $m=1$. In \cite{Tong2}, he showed that, except for $(D,p,x,m,n)=(2,3,7,3,1)$, if $(x,m,n)$ is a solution of \eqref{eq.2.19} with $m>1$, then either $(i)$ $D\equiv 1 \pmod{2}$, $p\equiv 1 \pmod{8}$, $m=2$ and $n=1$ or $(ii)$ $D\equiv 0 \pmod{2}$, $D$ has a prime divisor $d$ with $d\equiv 1 \pmod{2q}$, where $q$ is an odd prime divisor of $m$. However, Y.-F. He \cite{YFHe1} gave a counterexample that $(D,p,x,m,n)=(2,5,7,3,1)$, and found all solutions $(x,m,n)$ of \eqref{eq.2.19} with $m>1$. In addition, K. Gueth and L. Szalay \cite{GSz} solved the further equation
\begin{equation*}
x^2=3^2\pm 3\cdot 2^m+2^n, \, x,m,n\in\mathbb{N}.
\end{equation*} 

For fixed positive integers $D,k$ such that $\min\{D,k\}>1$ and $\gcd(2D,k)=1$, the equation
\begin{equation}\label{eq.2.20}
x^2+D^m=k^n, \, x,m,n\in\mathbb{N}
\end{equation}
is a naturally generalization of \eqref{eq.2.15}. So far there have very few results related to \eqref{eq.2.20}. In this respect, Y.-Z. Hu \cite{Hu1} proved that if $(D,k)=(3s^2-1,4s^2-1)$, where $s$ is a positive integer and $3s^2-1$ is an odd prime power, then \eqref{eq.2.20} has only the solutions $(x,m,n)=(s,1,1)$ and $(8s^3-3s,1,3)$. Y-Z. Hu and R.-X. Liu \cite{HuLiu} gave a similar result for $(D,k)=(3s^2+1,4s^2+1)$. J.-P. Wang and T.-T. Wang \cite{WangTang} proved that if $(D,k)=(3s^2+1,4s^2+1)$, where $s$ is a positive integer with $s\equiv 3 \pmod{6}$, then \eqref{eq.2.20} has only the solutions $(x,m,n)=(s,1,1)$ and $(8s^3+3s,1,3)$.

We now introduce two conjectures related to \eqref{eq.2.20}. Let $(a,b,c)$ be a primitive Pythagorean triple such that
\begin{equation}\label{eq.2.21}
a^2+b^2=c^2, \, a,b,c\in\mathbb{N},\, \gcd(a,b)=1,\, a\equiv 0\pmod{2}.
\end{equation}
In 1993, N. Terai \cite{Ter1} proposed the following conjecture.
\begin{conj}\label{conj.2.5}
The equation
\begin{equation}\label{eq.2.22}
x^2+b^y=c^z, \, x,y,z\in\mathbb{N}
\end{equation}
has only the solution $(x,y,z)=(a,2,2)$.
\end{conj}

Obviously, \eqref{eq.2.22} is the special case of \eqref{eq.2.20} for $(D,k)=(b,c)$. Conjecture \ref{conj.2.5} has been verified in the following cases:

$(i)$ (N. Terai \cite{Ter1}) $b\equiv 1 \pmod{4}$, $b^2+1=2c$, $b$ and $c$ are odd primes satisfying some conditions.

$(ii)$ (X.-G. Chen and M.-H. Le \cite{CL2}) $b\not\equiv 1 \pmod{16}$, $b^2+1=2c$, $b$ and $c$ are odd primes.

$(iii)$ (K.-Y. Chen \cite{KCh2}) $b$ is an odd prime with $b=16(3s+1)+1$, where $s$ is a positive integer.

$(iv)$ (M.-H. Le \cite{Le28}) $b\equiv \pm 3 \pmod{8}$, $b>8\cdot 10^6$, $c$ is an odd prime power.

$(v)$ (P.-Z. Yuan \cite{Yu2}; P.-Z. Yuan and J.-B. Wang \cite{YuWang}) $b\equiv \pm 3 \pmod{8}$, $c$ is an odd prime power.

$(vi)$ (M.-H. Le \cite{Le46}; L.-C. Gu \cite{Gu}; J.-Y. Hu and H. Zhang \cite{HuZh}) $b\equiv 7 \pmod{8}$, either $b$ or $c$ is an odd prime power.

$(vii)$ (S.-C. Yang \cite{SYang1}) $c$ is an odd prime power, $y\equiv z\equiv 0 \pmod{2}.$

In general, Conjecture \ref{conj.2.5} is not solved as yet. 

In 2001, as an analogue of Conjecture \ref{conj.2.5}, Z.-F. Cao and X.-L. Dong \cite{CaDo3} conjectured that if $a,b,c,l,m,n$ are fixed positive integers such that
\begin{equation}\label{eq.2.23}
a^l+b^m=c^n,\, \min\{a,b,c,l,m,n\}>1,\, \gcd(a,b)=1,\,a\equiv 0\pmod{2}
\end{equation}
then the equation
\begin{equation}\label{eq.2.24}
x^l+b^y=c^z, \, x,y,z\in\mathbb{N}
\end{equation}
has only the solution $(x,y,z)=(a,m,n)$. However, it has been found that there are infinitely many counterexamples to the above conjecture (see \cite{CDL}). Therefore, the condition $\min\{y,z\}>1$ is necessary in this conjecture, namely, the conjecture should be changed to the following:
\begin{conj}\label{conj.2.6}
If $a,b,c,l,m,n$ satisfy \eqref{eq.2.23}, then \eqref{eq.2.24} has only the solution $(x,y,z)=(a,m,n)$ with $\min\{y,z\}>1$.
\end{conj}      

Since Conjecture \ref{conj.2.6} is more involved than Conjecture \ref{conj.2.5}, the problem is far from solved. Except for the results on Conjecture \ref{conj.2.5} mentioned above, the most existing results concerning Conjecture \ref{conj.2.6} are limited to the case that $l=m=2$, $n\equiv 1 \pmod{2}$ and
\begin{equation}\label{eq.2.25}
a=s\left\lvert \sum\limits_{i=0}^{(n-1)/2}(-1)^{i}\binom {n} {2i}s^{n-2i-1}\right\rvert,\, b=\left\lvert \sum\limits_{i=0}^{(n-1)/2}(-1)^{i}\binom {n} {2i+1}s^{n-2i-1}\right\rvert, \, \\ c=s^2+1,\, s\in\mathbb{N},\, s\equiv 0\pmod{2}.
\end{equation}
In this respect, Conjecture \ref{conj.2.6} has been verified in the following cases:

$(i)$ (Z.-F. Cao, X.-L. Dong and Z.-Li \cite{CDL}) $b$ is an odd prime power with $b\equiv 3 \pmod{4}$.

$(ii)$ (M.-H. Le \cite{Le53}) $b\equiv 3\pmod{4}$, either $b$ or $c$ is an odd prime power.

$(iii)$ (W.-J. Guan \cite{Guan1}) $n=5$, $s$ is a power of 2, $y\equiv 0\pmod{2}$.

$(iv)$ (M.-Y. Lin \cite{Lin1}) $n\equiv 3\pmod{4}$, $s$ is a power of 2 with $s>n/\pi.$

$(v)$ (Y.-Z. Hu \cite{Hu4}) $n=5$, $b$ is an odd prime power.

$(vi)$ (Y.-Z. Hu and M.-H. Le \cite{HuLe1}) $n\equiv 5\pmod{8}$, $b$ is an odd prime power, $c>10^{12}n^4$.

In 2014, N. Terai \cite{Ter2} proposed another conjecture concerning \eqref{eq.2.22} as follows:
\begin{conj}\label{conj.2.7}
If $b,c,r,t$ are fixed positive integers such that
\begin{equation}\label{eq.2.26}
b^r+1=2c^t, \, \min\{b,c\}>1,\, t\in\{1,2\},
\end{equation}
then \eqref{eq.2.22} has only the solution $(x,y,z)=(c^t-1,r,2t)$.	
\end{conj}

This conjecture has been verified in the following cases:

$(i)$ (N. Terai \cite{Ter2}) $r=t=1$, $b$ is an odd prime, $c\le 30$ and $c\neq 12,24.$

$(ii)$ (R.-Q. Fu and H. Yang \cite{FuYa}) $r=t=1$, either $2c-1$ has a divisor $d$ with $d\equiv \pm 3 \pmod{8}$ or $ord_2(c-1)\equiv 1 \pmod{2}$, $z\equiv 1\pmod{2}.$ 

$(iii)$ (M.-J. Deng, J. Guo and A.-J. Xu \cite{DGX}) $r=t=1$, $3\le c\le 499$ and $c\equiv 3 \pmod{4}$.

$(iv)$ (S. Cenberci and H. Senay \cite{Cese}) $r=t=2$, $b,c$ are odd primes and $c$ satisfies some conditions.

$(v)$ (Y.-E. Zhao and X.-Q. Zhao \cite{ZZhao}) $r=t=2$, $b,c$ are odd primes with $b\equiv 3 \pmod{4}$.

$(vi)$ (G.-R. He \cite{He2}) $r=t=2$, $c\equiv \pm 3 \pmod{8}$.

$(vii)$ (X. Zhang \cite{Zh}) $r=t=2$, $b\equiv \pm 3 \pmod{8}$, $c$ is an odd prime.

$(viii)$ (J.-Y. Hu and X.-X. Li \cite{HuLi2}) $r>2$, $t=2$.

Finally, we introduce two problems related to \eqref{eq.2.20} in combinatorial mathematics. For odd primes $p$, S.-L. Ma \cite{LMa} proposed the following two conjectures in 1991.
\begin{conj}\label{conj.2.8}
The equation
\begin{equation}\label{eq.2.27}
x^2=2^{2a+2}p^{2n}-2^{2a+2}p^{m+n}+1, \, a\in\mathbb{Z}, \, a\ge 0, \, x,m,n\in\mathbb{N}
\end{equation}
has no solutions $(x,a,m,n)$.
\end{conj}
\begin{conj}\label{conj.2.9}
The equation
\begin{equation}\label{eq.2.28}
x^2=2^{2a+2}p^{2n}-2^{a+2}p^{m+n}+1, \, a\in\mathbb{Z}, \, a\ge 0, \, x,m,n\in\mathbb{N}
\end{equation}
has only the solution $(p,x,a,m,n)=(5,49,3,2,1)$.
\end{conj}

Five years later, Conjecture \ref{conj.2.8} was solved by M.-H. Le and Q. Xiang \cite{LeXi}. However, Conjecture \ref{conj.2.9} has not yet been solved. In this respect, F. Luca and P. St\v{a}nic\v{a} \cite{LucSt} proved that, for any fixed $p$, \eqref{eq.2.28} has at most $2^{30000}$ solutions $(x,a,m,n)$. M.-H. Zhu and T. Cheng \cite{ZCh} proved that \eqref{eq.2.28} has no solutions $(x,a,m,n)$ with $n\ge m$. Y.-F. He and Q. Tian \cite{YFHe2} further proved that the solutions $(x,a,m,n)$ of \eqref{eq.2.28} satisfy $m\ge 2n$ and $x=2^{a+1}s+\delta=2p^{2n}t-\delta$, where $\delta=(-1)^{(x-1)/2}$, $s,t$ are positive integers with $2^a-p^{m-n}=st.$ Moreover, \eqref{eq.2.28} has only the solution $(p,x,a,m,n)=(5,49,3,2,1)$ for $t=1$.

In addition, X.-L. Dong and Z.-F. Cao \cite{DoCa3}, X.-G. Guan \cite{Guan2}, J.-G. Luo, A. Togb\'{e} and P.-Z. Yuan \cite{LTY} discussed some equations of origin similar to \eqref{eq.2.27} and \eqref{eq.2.28}.

\section{The generalized Lebesgue-Nagell equation}

\subsection{The equation $x^2\pm D=y^n$}
\quad

\quad

In this subsection we first consider the solution of \eqref{eq.1.3} for $D\not\equiv 7 \pmod{8}$. Since $n>2$, then the solutions $(x,y,n)$ of \eqref{eq.1.3} must be satisfy $y\equiv 1 \pmod{2}$, so the solution of it is relatively easy. In 1993, on the basis of summarizing the previous works, J.H.E. Cohn \cite{Cohn5} solved \eqref{eq.1.3} for 77 values of $D$ with $D\not\equiv 7 \pmod{8}$ in the range $1 \le D\le 100$. Ten years later, using the BHV theorem, J.H.E. Cohn \cite{Cohn6} simplified his proof in \cite{Cohn5}. Afterwards, M. Mignotte and B. M. M. de Weger \cite{MigDe} solved \eqref{eq.1.3} for $D\in\{74,86\}$. So far, \eqref{eq.1.3} has been solved for $D\not\equiv 7 \pmod{8}$ with  $1 \le D\le 100$. In addition, M.-H. Le \cite{Le42} proved that if $D=163$, then \eqref{eq.1.3} has no solutions $(x,y,n)$.

As we pointed out in Introduction, if  $D\equiv 7 \pmod{8}$, then the solution of \eqref{eq.1.3} is very difficult. In this respect, using various methods including modular approach, S. Siksek and J. E. Cremona \cite{SiCr}, M.A. Bennett and C. M. Skinner \cite{BS} solved \eqref{eq.1.3} for $D=7$ and $D\in\{55,95\}$, respectively. In 2006, Y. Bugeaud, M. Mignotte and S. Siksek \cite{BMS2} solved \eqref{eq.1.3} for the remaining 19 values of $D$ in the range $1\le D \le 100$. In addition, for $D \le 10^{12}$, P.G. Walsh \cite{Wa2} completely found all solutions $(x,y,n)$ of \eqref{eq.1.3} with $n\equiv 0 \pmod{2}$.

For a general value of $D$, the solution of \eqref{eq.1.3} is a time-consuming and trivial work, so we consider about \eqref{eq.1.3} in a certain way.

For $D=a^2$, where $a$ is a positive integer with $a>1$, \eqref{eq.1.3} can be written as
\begin{equation}\label{eq.3.1}
x^2+a^2=y^2, \, x,y,n\in\mathbb{N},\, \gcd(x,y)=1,\, n>2.
\end{equation}
Since $x^2+a^2$ can be split in $\mathbb{Q}(\sqrt{-1})$, \eqref{eq.3.1} is one of the easier types to solve for \eqref{eq.1.3}. In this respect, M.-H. Le \cite{Le51} gave some formulas for all solutions $(x,y,n)$ of \eqref{eq.3.1} if $a$ is an odd prime, and completely determined all solutions for odd primes $a$ with $a<100$. Sz. Tengely \cite{Teng1} solved \eqref{eq.3.1} for $a \le 501$. M. Liang \cite{Liang2} proved that if $a$ is an odd prime with $a>7$, then \eqref{eq.3.1} has at most $\log a$ solutions $(x,y,n)$. For a general $a$, the upper bound on the solutions $(x,y,n)$ of \eqref{eq.3.1} has the following results:

$(i)$ (M.-Y. Lin \cite{Lin3}) $n<\max\{4\cdot 10^8, 4(\log a)/\log5\}.$

$(ii)$ (Sz. Tengely \cite{Teng1}) If $y>50000$, then every prime divisor $p$ of $n$ satisfies $p<\max\{9511, 4(\log a)/\log 50000\}$.

$(iii)$ (Sz. Tengely \cite{Teng2}) $p<\max\{4651, 1.85\log a\}$ where $p$ is a prime divisor of $n$.

$(iv)$ (X.-W. Pan \cite{Pan2}) For any positive number $\epsilon$, if $n>C_2(\epsilon)$, then $n<(2+\epsilon)(\log a)/\log y$.

Similarly, one has solved the equation
\begin{equation}\label{eq.3.2}
x^2+a^2=2y^n,\, x,y,n\in\mathbb{N}, \, \gcd(x,y)=1,\, n>2
\end{equation}
for the following cases:

$(i)$ (C. St\"{o}rmer \cite{St}) $a=1$.

$(ii)$ (I. Pink and Sz. Tengely \cite{PiSz}) $a\le 1000$ and $n\le 80$.

$(iii)$ (Sz. Tengely \cite{Teng1}) $a\le 305$.

$(iv)$ (L.-Y. Yang, J.-H. Chen and J.-L. Sun \cite{LYang}) $a\le 10^8.$

In general, I. Pink \cite{Pin1} proved that all solutions $(x,y,n)$ of \eqref{eq.3.2} satisfy $n<C_3(l,p)$, where $l,p$ are the number and the maximum of distinct prime divisors of $a$.

For a fixed positive integer $D$,
\begin{equation}\label{eq.3.3}
Dx^2+1=y^n, \, x,y,n\in\mathbb{N}, \, n>2
\end{equation}
is a class of equations that can easily be confused with \eqref{eq.1.3}. For example, M.G. Leu and G.W. Li \cite{LeLi} showed that if $D=2$, then \eqref{eq.3.3} has only the solutions $(x,y,n)=(2,3,3)$ and $(11,3,5)$ with $y=3$. Although the above result has been proved by W. Ljunggren \cite{Lj2} sixth years ago, but the last one of above solutions is often left out in some important literatures. In this respect, \eqref{eq.3.3} has been solved for the following cases:

$(i)$ (Z.-F. Cao \cite{Cao1}) $D=7.$

$(ii)$ (Q. Sun \cite{Sun}) $D\in\{15,23\}$.

$(iii)$ (Q.-C. Xiao \cite{Xia}) $D\in\{7,15,39,87\}$.

$(iv)$ (J.H.E. Cohn \cite{Cohn8}; E. Herrman, I. J\'{a}r\'{a}si and A. Peth\H{o} \cite{HJP}) $1\le D \le 100.$

In addition, P.-Z. Yuan \cite{Yu6} proved that if $(D,y)\neq (2,3), (6,7), (7,2)$ or $((t^2-2)/s^2,t^2-1)$, where $s,t$ are positive integers, then \eqref{eq.3.3} has at most one solution $(x,y,n)$ with fixed $y$.
Z.-F. Cao \cite{Cao2} and M.-H. Le \cite{Le10} discussed the solvability of \eqref{eq.3.3} for $n\equiv 2 \pmod{4}$. M.-H. Le \cite{Le35} gave an upper bound for the number of solutions of \eqref{eq.3.3}.

In 1981, Z. Ke and Q. Sun \cite{KeSun1,KeSun2} proved that if $D>2$ and $D$ has no prime divisors $p$ with $p\equiv 1 \pmod{6}$, then \eqref{eq.3.3} has no solutions $(x,y,n)$ with $n=3$. Ten years later, J.H.E. Cohn \cite{Cohn2} independently proved the same result. In 1998, M.-H. Le \cite{Le6} proved a similar result for $n=5$. The above results also hold true for the equation
\begin{equation}\label{eq.3.4}
Dx^2-1=y^n, \, x,y,n\in\mathbb{N}, \, n>2.
\end{equation}
In this respect, Z.-F. Cao \cite{Cao7} and J.H.E. Cohn \cite{Cohn10} proved that if $D=n=p$, where $p$ is an odd prime with $p\equiv 1 \pmod{4}$, and $B_{(p-1)/2}\not\equiv0 \pmod{p}$, where $B_{(p-1)/2}$ is the $\frac{1}{2}(p-1)-th$ Bernoulli number, then \eqref{eq.3.4} has no solutions $(x,y,n)$. M.A. Bennett and C.M. Skinner \cite{BS} proved that if $D=2$, then \eqref{eq.3.4} has only the solution $(x,y,n)=(78,23,3)$. In addition, M.A. Bennett \cite{Ben1} found all solutions $(x,y,n)$ of the equation
\begin{equation*}
x^2-1=Dy^{2n}, \, x,y,n\in\mathbb{N}, n>2
\end{equation*}   
for $1<D<100$.

Let $a,b$ be fixed positive integers. Obviously, if $(z,n)$ is a solution of the equation
\begin{equation}\label{eq.3.5}
(a^n-1)(b^n-1)=z^2, \, z,n\in\mathbb{N}, \, n>2,
\end{equation}
then we have
\begin{equation}\label{eq.3.6}
a^n-1=Dx^2,\, b^n-1=Dy^2, \, z=Dxy,\, D,x,y\in\mathbb{N}.
\end{equation}
We see from \eqref{eq.3.6} that the solution of \eqref{eq.3.5} is directly related to \eqref{eq.3.3}. In this respect, we have the following results:

$(i)$ (L. Szalay \cite{Sz}) If $(a,b)=(2,3)$ or $(2,5)$, then \eqref{eq.3.5} has no solutions $(z,n)$. If $(a,b)=(2,2^r)$, where $r$ is a positive integer with $r>1$, then \eqref{eq.3.5} has only the solution $(r,z,n)=(2,21,3)$.

$(ii)$ (L. Hajdu and L. Szalay \cite{HS}) If $(a,b)=(2,6)$ or $(a,a^r)$, where $a>2$ and $r>1$, then \eqref{eq.3.5} has no solutions $(z,n).$

$(iii)$ (J.H.E. Cohn \cite{Cohn7}) If $2\le a<b\le 12$, then \eqref{eq.3.5} is solved, except for 24 pairs of $(a,b)$.

$(iv)$ (F. Luca and P.G. Walsh \cite{LW}) If $2\le a<b\le 100$, then \eqref{eq.3.5} is solved, except for 69 pairs of $(a,b)$ included $(a,b)=(4,13)$ and $(13,28)$.

$(v)$ (Z.-J. Li and M. Tang \cite{LiTang1,LiTang2}) If $(a,b)=(4,13)$ or $(13,28)$, then \eqref{eq.3.5} has no solutions $(z,n).$

$(vi)$ (J.H.E. Cohn \cite{Cohn7}) \eqref{eq.3.5} has only the solution $(a,b,z,n)=(13,239,9653280,4)$ with $n\equiv \pmod{4}$.

$(vii)$ (M.-H. Le \cite{Le57}) If $a=2$ and $b\equiv 0 \pmod{3}$, then \eqref{eq.3.5} has no solutions $(z,n)$.

$(viii)$ (Z.-J. Li \cite{ZLi}) If $a=2$, $b\equiv 3, 19, 67, 83, 131, 147, 171$ or $179\pmod{200}$, then \eqref{eq.3.5} has no solutions $(z,n)$.

$(ix)$ (L. Li and L. Szalay \cite{LiSz}) If $a\equiv 2 \pmod{6}$ and $b\equiv 0 \pmod{3}$, then \eqref{eq.3.5} has no solutions $(z,n)$.

$(x)$ (M. Tang \cite{Tang}) If $a\equiv 0 \pmod{2}$ and $b\equiv 15 \pmod{20}$, then \eqref{eq.3.5} has no solutions $(z,n)$.

$(xi)$ (S.-C. Yang, W.-Q. Wu and H. Zheng \cite{YWZ}) If $(a,b)\equiv (5,8) \pmod{13}$ or $(6,7) \pmod{17}$, and $a,b$ satisfy some conditions, then \eqref{eq.3.5} has no solutions $(z,n)$.

$(xii)$ (P.-Z. Yuan and Z.-F. Zhang \cite{YuZh}) If $(a,b)\equiv (2,0) \pmod{3}$, $(3,0) \pmod{4}$, $(3,2) \pmod{4}$ or $(4,0) \pmod{5}$, then \eqref{eq.3.5} has no solutions $(z,n)$.

$(xii)$ (K. Ishii \cite{Ish2}) If $a\equiv 5 \pmod{6}$, $b\equiv 0 \pmod{3}$, and $a,b$ satisfy some conditions, then \eqref{eq.3.5} has no solutions $(z,n)$.

$(xiv)$ (J.H.E. Cohn \cite{Cohn7}) If $ord_2(a-1)\not\equiv ord_2(b-1) \pmod{2}$, 
then \eqref{eq.3.5} has no solutions $(z,n)$ with $n\equiv 1 \pmod{2}$.

$(xv)$ (X.-Y. Guo \cite{Guo}) All solutions $(z,n)$ of \eqref{eq.3.5} are determined for $ord_2(a-1)\not\equiv ord_2(b-1) \pmod{2}$.

$(xvi)$ (M. Liang \cite{Liang1}) If $a\equiv 2$ or $3 \pmod{4}$, $b=a+1$, then \eqref{eq.3.5} has no solutions $(z,n)$.

$(xvii)$ (H. Yang, Y.-T. Pei and R.-Q. Fu \cite{YPF}) If $a\equiv 5$ or $9 \pmod{16}$, $b=a+1$, then \eqref{eq.3.5} has no solutions $(z,n)$.

In 2002, J.H.E. Cohn \cite{Cohn7} proposed the following conjecture:
\begin{conj}
All solutions $(z,n)$ of \eqref{eq.3.5} satisfy $n\le 4$.
\end{conj}

This problem has not yet been solved. In addition, E.G. Goedhart and H.G. Grundman \cite{GG2}, G.-R. He \cite{He1}, M.-H. Le \cite{Le36}, M. Tang \cite{Tang}, P.G. Walsh \cite{Wa1}, Z.-F. Zhang \cite{ZZh} discussed some variations of \eqref{eq.3.5}, respectively.

Finally, we should point out that, for a fixed positive integer $D$ which is not a square, the solution of equations with the form
\begin{equation}\label{eq.3.7}
x^2-D=y^n,\, x,y,n\in\mathbb{N}, \, \gcd(x,y)=1, \, n>2
\end{equation}
is very difficult. For example, one only know the upper bound for solutions $(x,y,n)$ of the equation
\begin{equation}\label{eq.3.8}
x^2-2=y^n,\, x,y,n\in\mathbb{N},\, \gcd(x,y)=1,\, n>2
\end{equation}
(see \cite{Bug2}, \cite{GL} and \cite{Si2}), but cannot determine all of its solutions. Then we have
\begin{conj}\label{conj.3.2}
\eqref{eq.3.8}  has no solutions $(x,y,n)$.
\end{conj}

In addition, C.F. Barros \cite{Bar} gave some upper bounds for solutions $(x,y,n)$ of \eqref{eq.3.7} in the range $1<D<100$.

For $\delta\in\{1,2,4\}$, let $D_1,D_2$ be fixed integers such that $\gcd(D_1,D_2)=\gcd(\delta,D_1D_2)=1$ and $D_1D_2$ is not a square. Then the equation
\begin{equation}\label{eq.3.9}
D_1x^2+D_2=\delta y^n, \, x,y,n\in\mathbb{N},\, \gcd(x,y)=1,\, n>2
\end{equation}
is a naturally generalization of \eqref{eq.1.3}. By the way, the research of \eqref{eq.3.9} mainly focuses on discussing its solutions suitable some conditions concerning the divisibility of class number of the quadratic field $\mathbb{Q}(\sqrt{-D_1D_2})$. See subsection 3.3 of this survey for the results.      

\subsection{The equation $x^2\pm D^m=y^n$}
\quad

\quad

For fixed positive integers $D$ with $D>1$, the equation
\begin{equation}\label{eq.3.10}
x^2+D^m=y^n, \, x,y,m,n\in\mathbb{N},\, \gcd(x,y)=1, \, n>2
\end{equation}
is an exponential generalization of \eqref{eq.1.3}. This equation has been solved for some smaller values of $D$. For example, J.H.E. Cohn \cite{Cohn3,Cohn4} proved that if $D=2$, then \eqref{eq.3.10} has only the solutions $(x,y,m,n)=(5,3,1,3)$ and $(7,3,5,4)$ with $m\equiv 1 \pmod{2}$. Simultaneously, he conjectured that if $D=2$, then \eqref{eq.3.10} has only the solution $(x,y,m,n)=(11,5,2,3)$ with $m\equiv 0 \pmod{2}$. Around 1997, using the Baker method, M.-H. Le \cite{Le37,Le48} confirmed Cohn's conjecture. 
Afterwards, J.-Y. Hu and X.-X. Li \cite{HuLi1} used the BHV theorem to give a new proof of Le's result. In addition, \eqref{eq.3.10} has been solved in the following cases:

$(i)$ (S.A. Arif and F.S. Abu Muriefah \cite{ArAb1}; F. Luca \cite{Luc1}; L.-Q. Tao \cite{Tao1}; Y.-Z. Hu \cite{Hu3}) $D=3$, $(x,y,m,n)=(10,7,5,3)$ and $(46,13,4,3)$.

$(ii)$ (F.S. Abu Muriefah \cite{Abu4}; F.S. Abu Muriefah and S.A. Arif \cite{Abu6}; L.-Q. Tao \cite{Tao2}) $D=5$, no solutions.

$(iii)$ (F. Luca and A. Togb\'{e} \cite{LucTog}) $D=7$, $m\equiv 0 \pmod{2}$, $(x,y,m,n)=(24,5,2,4)$ and $(524,65,2,3)$.

$(iv)$ (I.N. Cangul, G. Soydan and Y. Simsek \cite{CSS}; G. Soydan, M. Demirci and I. N. Cangul \cite{SDC}) $D=11$, $(x,y,m,n)=(2,5,2,3),(4,3,1,3)$, $(58,15,1,3)$ and $(9324,443,3,3)$.

$(v)$ (B. Peker and S. Cenberci \cite{PekCen}) $D=19$, $(x,y,m,n)=(18,7,1,3)$ and $(22434,55,1,5)$.

$(vi)$ (N. Saradha and A. Srinivasan \cite{SarSr1}; H.-L. Zhu and M.-H. Le \cite{ZhuLe}) $D=67,$ $(x,y,m,n)=(110,23,1,3)$. $D\in\{43,167\}$, no solutions.

For $D=p^2$, where $p$ is an odd prime, \eqref{eq.3.10} can be written as
\begin{equation}\label{eq.3.11}
x^2+p^{2m}=y^n, \, x,y,m,n\in\mathbb{N}, \, \gcd(x,y)=1, \, n>2.
\end{equation}   
In this respect, M.-H. Le \cite{Le47} determined all solutions $(x,y,m,n)$ of \eqref{eq.3.11} with $n\equiv 0 \pmod{2}$. A. B\'{e}rczes and I. Pink \cite{BP1} solved \eqref{eq.3.11} for $p<100$. S.A. Arif and F.S. Abu Muriefah \cite{ArAb2} gave some solvability conditions of \eqref{eq.3.11}. X.-W. Pan \cite{Pan1} described all solutions $(x,y,m,n)$ of \eqref{eq.3.11}.

Similarly, Sz. Tengely \cite{Teng3} discussed the equation
\begin{equation}\label{eq.3.12}
x^2+p^{2m}=2y^n, \, x,y,m,n\in\mathbb{N}, \, \gcd(x,y)=1,\, n>2.
\end{equation}
He used the Baker theory to prove that there exist only finitely many odd primes $p$ such that \eqref{eq.3.12} has the solutions $(x,y,m,n)$ satisfying the following condition:
\begin{equation}\label{eq.3.13}
\textrm{m is an odd prime with} \, n>3, \, y \,\, \textrm{is not a sum of two consecutive squares}. 
\end{equation}
Five years later, using the BHV theorem, H.-L. Zhu, M.-H. Le and A. Togb\'{e} \cite{ZLT} desribed all solutions of \eqref{eq.3.12}. In the same paper, they proved that, for any odd prime $p$, \eqref{eq.3.12} has no solutions $(x,y,m,n)$ with \eqref{eq.3.13}.

In 1978, S. Rabinowitz \cite{Rab2} proved that the equation
\begin{equation}\label{eq.3.14}
nx^2+2^m=y^n,\, x,y,m,n\in\mathbb{N},\, \gcd(x,y)=1,\, n>2
\end{equation}
has only the solution $(x,y,m,n)=(21,11,3,3)$ with $n=3$. Afterwards, M.-H. Le \cite{Le32}, H.-M. Wu \cite{Wu1}, Y.-X. Wang and T.-T. Wang \cite{WW}, G. Soydan and I. N. Cangul \cite{SC} discussed various cases of \eqref{eq.3.14}, respectively. This equation was finally solved by F. Luca and G. Soydan \cite{LucSoy}, they proved that \eqref{eq.3.14} has only the solution $(x,y,m,n)=(21,11,3,3)$. In addition, E.G. Goedhart and H.G. Grundman \cite{GG1} proved that the equation
\begin{equation*}
nx^2+2^{m_1}3^{m_2}=y^n,\, x,y,m_1,m_2,n\in\mathbb{N},\, \gcd(x,y)=1,\, n>2
\end{equation*}  
has no solutions $(x,y,m_1,m_2,n).$

Let $p,q$ be fixed distinct odd primes. In 2000, F.S. Abu Muriefah \cite{Abu1} proved that if $p\not\equiv 7\pmod{8}$ and $q=3$, then the equation  
\begin{equation}\label{eq.3.15}
px^2+q^{2m}=y^p,\, x,y,m\in\mathbb{N},\, \gcd(x,y)=1
\end{equation}
has no solutions $(x,y,m)$. Eight years later, F.S. Abu Muriefah \cite{Abu5} showed that the above conclusion is true for any $q$. However, D. Goss \cite{Goss} pointed out that the proof of \cite{Abu5} is wrong. Afterwards, A. Larajdi, M. Mignotte and N. Tzanakis \cite{LMT} gave a correct proof of the above conclusion for $p\equiv 3 \pmod{4}$. X.-Y. Wang and H. Zhang \cite{WZ} solved the case that $p\equiv 1 \pmod{4}$, either $q<4p-1$ or $q<149$. But, so far it is known whether the above conclusion is correct for any fixed $p$ with $p\equiv 1 \pmod{p}.$ 

For fixed distinct primes $p_1,...,p_t$ $(t>1)$, let $S(p_1,\cdots,p_t)=\{p_1^{r_1}\cdots p_t^{r_t}|\, r_i\in\mathbb{Z}, \, r_i\ge 0,\, i=1,\cdots,t$, $(r_1,\cdots, r_t)\neq (0,\cdots, 0)\}$.
Then the equation 
\begin{equation}\label{eq.3.16}
x^2+D=y^n,\, D\in S(p_1,\cdots,p_t),\, x,y,n\in\mathbb{N},\, \gcd(x,y)=1,\, n>2
\end{equation}
is another exponential generalization of \eqref{eq.1.3}. This equation has been solved in the following cases:

$(i)$ (F. Luca \cite{Luc2}) $D\in S(2,3).$

$(ii)$ (F. Luca and A. Togb\'{e} \cite{LT1}) $D\in S(2,5).$

$(iii)$ (I.N. Cangul, M. Demirci, F. Luca, \'{A}. Pint\'{e}r and G. Soydan \cite{CDLPS}) $D\in S(2,11).$

$(iv)$ (F. Luca and A. Togb\'{e} \cite{LT2}) $D\in S(2,13).$

$(v)$ (S. Gou and T.-T. Wang \cite{GT}; A. Dabrowski \cite{Dab}) $D\in S(2,17).$

$(vi)$ (G. Soydan, M. Ulas and H.-L. Zhu \cite{SUZ}) $D\in S(2,19).$

$(vii)$ (A. Dabrowski \cite{Dab}) $D\in S(2,29)$ or $S(2,41).$

$(viii)$ (H.-L. Zhu, M.-H. Le and G. Soydan \cite{ZLS}; H.-L. Zhu, G. Soydan and Q. Wei \cite{ZSW}) $D\in S(2,p)$ where $p$ is an odd prime, $n=4$.

$(ix)$ (H.-L. Zhu, M.-H. Le, G. Soydan and A. Togb\'{e} \cite{ZLST}) $D\in S(2,p)$ where $p$ is an odd prime satisfying some conditions.

$(x)$ (Y.-G. Ma \cite{YMa}) $D\in S(3,11)$ and $n\in\{3,5\}.$

$(xi)$ (I.N. Cangul, M. Demirci, G. Soydan and N. Tzanakis \cite{CDST}; G. Soydan and N. Tzanakis \cite{ST}) $D\in S(5,11).$

$(xii)$ (F.S. Abu Muriefah, F. Luca and A. Togb\'{e} \cite{Abu10}) $D\in S(5,13).$

$(xiii)$ (I. Pink and Z. R\'{a}bai \cite{PRab}; Y.-H. Gao \cite{Gao}) $D\in S(5,17).$

$(xiv)$ (G. Soydan \cite{Soy}) $D\in S(7,11).$

$(xv)$ (A. B\'{e}rczes and I. Pink\cite{BP3}) $D\in S(11,17).$

$(xvi)$ (I.N. Cangul, M. Demirci, I. Inam, F. Luca and G. Soydan \cite{CDILS}) $D\in S(2,3,11).$

$(xvii)$ (H. Godinho, D. Marques and A. Togb\'{e} \cite{GMT2}) $D\in S(2,3,17)$ or $S(2,13,17).$ 

$(xviii)$ (E. Goins, F. Luca and A. Togb\'{e} \cite{GLT}) $D\in S(2,5,13).$

$(xix)$ (H. Godinho, D. Marques and A. Togb\'{e} \cite{GMT1}) $D\in S(2,5,17).$

$(xx)$ (I. Pink \cite{Pin2}) $D\in\ S(2,3,5,7)$ and $y\equiv 1 \pmod{2}$ if $D=7s^2$ or $15s^2$, where $s$ is a positive integer.

In general, M.-Y. Lin \cite{Lin2} proved that all solutions $(x,y,n)$ of \eqref{eq.3.16} satisfy $n<\dfrac{120}{\pi}p_1\cdots p_t \log(2ep_1\cdots p_t)$.

Compared with \eqref{eq.3.10}, the solution of its daul form
\begin{equation}\label{eq.3.17}
x^2-D^m=y^n,\, x,y,m,n\in\mathbb{N},\, \gcd(x,y)=1,\, n>2
\end{equation}
is a more difficult problem. For example, in 1977, S. Rabinowitz \cite{Rab1} used the algebraic number theory method to proved that the equation
\begin{equation}\label{eq.3.18}
x^2-2^m=y^n,\, x,y,m,n\in\mathbb{N},\, \gcd(x,y)=1,\, n>2
\end{equation} 
has only the solution $(x,y,m,n)=(71,17,7,3)$ with $n=3.$ Around 1996, using the Baker method, M.-H. Le \cite{Le34}, Y.-D. Guo and M.-H. Le \cite{GL}, Y. Bugeaud \cite{Bug2} proved that all solutions $(x,y,m,n)$ of \eqref{eq.3.18} satisfy $\max\{x,y,m,n\}<C_5$. Until 2004, using modular approach, M.A. Bennett and C.M. Skinner \cite{BS} proved that \eqref{eq.3.18} has only the solution $(x,y,m,n)=(71,17,7,3)$ with $m>1$. In addition, F.S. Abu Muriefah and A. Al-Rashed \cite{Abu7}, M. Bauer and M.A. Bennett \cite{BB2}, Y. Bugeaud \cite{Bug1}, Y.-Y. Liu \cite{YLiu}, Z.-W. Liu \cite{ZLiu2} discussed some other cases of \eqref{eq.3.17}, respectively.

Finally, we introduce two polynomial-exponential equations similar to \eqref{eq.3.10} and \eqref{eq.3.17}. In 2000, F. Luca and M. Mignotte \cite{LuMig} used the Baker method to prove that the equation
\begin{equation}\label{eq.3.19}
x^y+y^x=z^2,\, x,y,z\in\mathbb{N},\, \min\{x,y\}>1, \,\gcd(x,y)=1
\end{equation} 
has no solutions $(x,y,z)$ with $xy\equiv 0 \pmod{2}$, and the equation
\begin{equation}\label{eq.3.20}
x^y-y^x=z^2,\, x,y,z\in\mathbb{N},\, \min\{x,y\}>1, \,\gcd(x,y)=1
\end{equation} 
has only the solution $(x,y,z)=(3,2,1)$ with $xy\equiv 0 \pmod{2}$. Seven years later, using the BHV theorem, M.-H. Le \cite{Le54} completely solved \eqref{eq.3.20}. He proved that \eqref{eq.3.20} has no solutions $(x,y,z)$ with $xy\equiv 1 \pmod{2}$. Afterwards, X.-Y. Wang \cite{Wang}, H. Yang and R.-Q. Fu \cite{YanFu2} gave some necessary conditions for \eqref{eq.3.19} to have solutions $(x,y,z)$ with $xy\equiv 1 \pmod{2}$.

In 2012, Z.-F. Zhang, J.-G. Luo and P.-Z. Yuan \cite{ZLY2} discussed a generalization of \eqref{eq.3.19} with the form  
\begin{equation}\label{eq.3.21}
x^y+y^x=z^z,\, x,y,z\in\mathbb{N}.
\end{equation}
Using the baker method, they proved that all solutions $(x,y,z)$ of \eqref{eq.3.21} satisfy $z<2.8\times10^9$. In the same paper, the authors conjectured that \eqref{eq.3.21} has no solutions $(x,y,z)$ with $\min\{x,y,z\}>1$. In this respect, Y. Deng and W.-P. Zhang \cite{DZW}, X.-Y. Du \cite{Du}, H.-M. Wu \cite{Wu2} proved that \eqref{eq.3.21} has no solutions for some cases. Very recently, M. Cipu \cite{Cipu} completely confirmed the conjecture of \cite{ZLY2}. In addition, Z.-F. Zhang, J.-G. Luo and P.-Z. Yuan \cite{ZYu} discussed other variations of \eqref{eq.3.19} and \eqref{eq.3.20}.        
\subsection{The divisibilty of class number of $\mathbb{Q}(\sqrt{\pm (a^2-\delta c^n)})$}
\quad

\quad

For any integer $d$ with square free, let $h(d)$ denote the class number of the quadratic field $\mathbb{Q}(\sqrt{d})$. There have many explicit results on the divisibility of $h(d)$  are directly related to the solution of \eqref{eq.1.3} and its generalizations. The research on this can be traced back to the works of T. Nagell \cite{Nag1} and W. Ljunggren \cite{Lj1}. In this respect, the early results can be seen in \cite{Le22,Le23,Le24} and \cite{Le31}. We now break down some recent result into the following two cases.

Case I. Let $a,b,c,n,\delta$ be fixed positive integers such that 
\begin{equation}\label{eq.3.22}
a^2+db^2=\delta c^n,\,  \gcd(a,db)=1,\, n>2,\, \delta\in\{1,4\},\, c>1,\,c\equiv 1\pmod{2} \, \textrm{if} \, \delta=1.
\end{equation}

$(i)$ (M. J. Cowles \cite{Cow}) If $a=b=1$, $c$ and $n$ are odd primes, $\delta=4$, then $h(-d)\equiv 0 \pmod{n}$.

$(ii)$ (B.H. Gross and D.E. Rohrlich \cite{GR}) If $a=1$, $n$ is an odd prime, $\delta=4$, then $h(-d)\equiv 0 \pmod{n}$.

$(iii)$ (H.-W. Lu \cite{Lu1}) If $a=b=1$, $\delta=4$, then $h(-d)\equiv 0 \pmod{n}$.

$(iv)$ (S.R. Louboutin \cite{Lou}) If $d>3$, $a=1$, $c\equiv 1 \pmod{2}$, $c$ has an odd prime divisor $p$ with $p\equiv 3 \pmod{4}$, $\delta=4$, then $h(-d)\equiv 0 \pmod{n}$.

$(v)$ (M.-H. Le \cite{Le4}) If $d>3$, $a=1$, $\delta=4$, then 
 
\begin{equation*}
h(-d)\equiv
\begin{cases}
0 \pmod{\frac{n}{4}}, \quad (d,b,c,n)=(7,3,2,4),\\
0 \pmod{\frac{n}{2}}, \quad n\equiv 0 \pmod{2},\; b=b_1b_2,\; \textrm{and}\\
\; \quad \quad \quad \quad \quad \quad \; \; \;       b_1^2-db_2^2=(-1)^c2 \,,b_1,b_2\in\mathbb{N},\\
0 \pmod{n}, \quad \textrm{otherwise}.
\end{cases}
\end{equation*}

$(vi)$ (X.-B. Zhang and X.-X. Li \cite{ZhLi}) If $d>3$, $a=1$, $\delta=4$, then a necessary and sufficent condition for $h(-d)\not\equiv 0 \pmod{n}$ is given.

$(vii)$ (N.C. Ankeny and S. Chowla \cite{AnCh}) If $a\equiv 0 \pmod{2}$, $c=3$, $a<\sqrt{2.3^{n-1}}$, $\delta=1$, then $h(-d)\equiv 0\pmod{n}$.

$(viii)$ (Y. Kishi \cite{Kish}) If $a=2^m$, where $m$ is a positive integer with $(m,n)\neq (2,3)$, $c=3$, $\delta=1$, then $h(-d)\equiv 0 \pmod{n}.$

$(ix)$ (A. Ito \cite{Ito1}) If $a=2^m$, where $m$ is a positive integer with $m\equiv 1 \pmod{2}$, $c$ is an odd prime, $\delta=1$, then $h(-d)\equiv 0 \pmod{n}.$ 

$(x)$ (M.-H. Zhu and T.-T. Wang \cite{ZWang}) If $a=2^m$, where $m$ is a positive integer, $\delta=1$, then $h(-d)\equiv 0 \pmod{n}$ or $\pmod{\dfrac{n}{3}}.$ 

$(xi)$ (A. Ito \cite{Ito3}) If $a=3^m$, where $m$ is a positive integer, $c$ is a prime, $n\not\equiv 2 \pmod{4}$, $\delta=4$, then $h(-d)\equiv 0 \pmod{n}.$

$(xii)$ (K. Ishii \cite{Ish1}; A. Ito \cite{Ito2}) If $d\equiv 1,5$ or $7 \pmod{8}$, $n\not\equiv 0 \pmod{3}$, then $h(-d)\equiv 0 \pmod{n}.$

$(xiii)$ (R.A. Mollin \cite{Mol1}) If $a<\sqrt{\delta(c-1)c^{n-1}},$ $b=1$, then $h(-d)\equiv 0 \pmod{n}$.

$(xiv)$ (M.-H. Le \cite{Le8}) If $d>exp(exp(exp(1000)))$, $b=1$, $\delta=4$, then  
\begin{equation*}
h(-d)\equiv
\begin{cases}
0 \pmod{\frac{n}{2}}, \quad n\equiv 0 \pmod{2},\; d=4c^{n/2}-1,\; a=2c^{n/2}-1,\\
0 \pmod{n}, \quad \textrm{otherwise},
\end{cases}
\end{equation*}
except for $(d,a,c,n)=(3(2s^2+1)^2\pm 4,2s+1,(2s+1)^2\pm 1,3)$ where $s$ is a positive integer.

In addition, F.S. Abu Muriefah \cite{Abu2,Abu3}, F.S. Abu Muriefah, F. Luca, S. Siksek and Sz. Tengely \cite{Abu9}, F. Luca, Sz. Tengely and A. Togb\'{e} \cite{LTT}, S.A. Arif and F.S. Abu Muriefah \cite{ArAb2,ArAb3}, S.A. Arif and A.S. Al-Ali \cite{ArAl1,ArAl2}, A. B\'{e}rczes and I. Pink \cite{BP2}, Y.F. Bilu \cite{Bil}, Y. Bugeaud \cite{Bug3}, Y. Bugeaud and T.N. Shorey \cite{BugSh}, Z.-F. Cao \cite{Cao3,Cao4,Cao5,Cao8}, Z.-F. Cao and X.-L. Dong \cite{CaDo1,CaDo2}, M. Cipu \cite{Cip}, J.H.E. Cohn \cite{Cohn9}, X.-L. Dong and Z.-F. Cao \cite{DoCa2}, M.-H. Le \cite{Le40,Le41,Le43}, M. Mignotte \cite{Mig}, R.A. Mollin \cite{Mol2,Mol3}, H.-L. Zhu  \cite{Zhu} discussed the solution of \eqref{eq.3.9} and its generalizations concerning the divisibility of class numbers of imaginary quadratic fields.

Case II. Let $a,c,n$ be fixed positive integers such that
\begin{equation}\label{eq.3.23}
da^2-1=4c^{2n},\, c>1,\, n>2.
\end{equation}

$(i)$ (P.J. Weinberger \cite{Wei}; R.A. Mollin and H.C. Williams \cite{MolWil}) For any fixed, there exist infinitely many values of $c$ which make $h(d)\equiv 0 \pmod{n}.$

$(ii)$ (H.-W. Lu \cite{Lu2}) If $a=1$, then $h(d)\equiv 0 \pmod{n}.$

$(iii)$ (M.-H. Le \cite{Le12}) If prime divisor $p$ of $n$ and $q$ of $c$ satisfy $\gcd(p,(q-1)q)=1$, $(U_1,V_1)=(2c^n,a)$ is the least solution of \eqref{eq.2.10} for $D=d$, then $h(d)\equiv 0 \pmod{n},$ except for $(d,a,c,n)=(41,5,2,4)$.

$(iv)$ (Z.-F. Cao \cite{Cao5}) If $a\le c^{n/2}$, $c\equiv 1 \pmod{2}$, $(U_1,V_1)=(2c,a)$ is the least solution of \eqref{eq.2.10} for $D=d$, then $h(d)\equiv 0 \pmod{n}.$

$(v)$ (P.-Z. Yuan \cite{Yu1}; X.-L. Dong and Z.-F. Cao \cite{DoCa1}) If prime divisors $p$ of $n$ and $q$ of $c$ satisfy $\gcd(p,q^2-1)=1$, then $h(d)\equiv 0 \pmod{n}.$

$(vi)$ (X.-L. Dong and Z.-F. Cao \cite{DoCa1}) If either $a\le \frac{1}{2}c^{0.4226n}$ or every prime divisor of $a$ divides $d$, then $h(d)\equiv 0 \pmod{n}.$

In addition, X.-K. Zhang \cite{XKZh2,XKZh3,ZhWas} proved some similar results concerning the divisibility of class numbers of real quadratic fields.

\bigskip
\subsection*{Acknowledgements}
The authors would like to thank Professor Mike Bennett and Dr. Paul Voutier for reading the original manuscript carefully and giving valuable advice. The second author was supported by T\"{U}B\.{I}TAK (the Scientific and Technological Research Council of Turkey) under Project No: 117F287.

\end{document}